\documentclass[11pt]{article}
\usepackage[utf8]{inputenc}
\usepackage[T1]{fontenc}
\usepackage{graphicx}
\usepackage{amsmath,amsthm,amssymb}
\usepackage{geometry}
\usepackage{url}
\usepackage{hyperref}
\usepackage{float}
\usepackage{longtable}
\usepackage{booktabs}
\usepackage{multirow}
\usepackage[numbers]{natbib}
\usepackage[x11names]{xcolor}
\usepackage{textcomp}
\usepackage{tikz}
\usepackage{caption}

\newtheorem{definition}{Definition}

\newtheorem{remark}{Remark}

\newtheorem{fact}{Fact}

\newtheorem{observation}{Observation}

\title{Algebraic Classification of All 880 Fourth-Order Magic Squares and the Discovery of Complete Alternating Magic Squares}
\author{Kenichi Takemura\\
\textit{Independent Researcher, Tokyo, Japan }\\
\textit{Email: mahora1123@gmail.com}\\
}
\date{\today}

\begin{document}

\maketitle

\begin{abstract}
In this paper, we introduce a newly defined algebraic invariant for square matrices termed the \emph{Alternating Power Difference (APD)}. The APD is defined as the signed sum of the powers of diagonal sums along permutations of the symmetric group, distinguishing between even and odd permutations. It serves as a measure of the broken even-odd symmetry inherent in a matrix through higher-order moments.

We applied this invariant to all 880 essentially different normal $4\times4$ magic squares (excluding symmetries) and defined the \emph{First Appearance Degree} $m_1$ as the minimum power at which the APD first becomes non-zero. Through an exhaustive computational search, we found that these magic squares are categorized into three clearly separated classes: $m_1=3$ (240 squares), $m_1=4$ (624 squares), and $m_1=\infty$ (16 squares).

In particular, the case $m_1=\infty$ identifies exceptionally rare magic squares for which the APD vanishes at all degrees. We refer to these as \emph{Complete Alternating Magic Squares} and demonstrate that they possess a strong algebraic symmetry undetectable by conventional geometric classifications or link-line patterns. Furthermore, we reveal that the APD-based classification refines the classical link-line classification based on complementary sum pairs, showing that each geometric type is clearly distinguished by its first appearance degree.

All results in this paper are based on exhaustive computations and are fully reproducible. Our findings suggest that the APD is an effective new invariant for detecting hidden algebraic structures in magic squares and related combinatorial matrices.
\end{abstract}

% --- Keywords ---
\providecommand{\keywords}[1]{\textbf{\textit{Keywords:}} #1}


% --- MSC Codes (Mathematics Subject Classification 2020) ---
% 05B15: Orthogonal arrays, Latin squares, Arithmetical magic squares
% 15B99: Special matrices (None of the above, but in this section)
% 00A08: Recreational mathematics
% 68W30: Symbolic computation and algebraic computation
\makeatletter
\newcommand{\msc}[1]{\textbf{\textit{MSC2020:}} #1}
\makeatother

\msc{05B15, 15B99, 00A08, 68W30}

\section{Definition and Total Number of Magic Squares}

\begin{definition}[Fourth-Order Magic Square]
A fourth-order magic square is an arrangement of the distinct integers from 1 to 16 in a $4\times4$ square grid such that the sum of the numbers in each row, each column, and both main diagonals is the same. This sum is called the magic constant, which is 34 for a fourth-order magic square.
\end{definition}

\begin{fact}
The total number of fourth-order magic squares is 880. This count represents essentially different arrangements, excluding those obtainable by rotation or reflection \cite{trump2018,andrews1976}.
\end{fact}

In this study, we classify all 880 types of magic squares using the concept of link-lines. The magic square data analyzed is provided in CSV format, with each square assigned an index from 1 to 880. To ensure research transparency and reproducibility, all classification results and analysis programs used for numerical verification are available in the author's GitHub repository:

\begin{quote}
    \url{https://github.com/soaisu-ken/magic-square-apd-classification/releases/tag/v1.0.0}
\end{quote}

All major numerical results presented in this paper can be independently reproduced using the code in the aforementioned repository.

\section{Concept of Link-Lines}

\subsection{Definition of Link-Lines}

\begin{definition}[Link-Line]
In a fourth-order magic square, a pair of numbers whose sum is 17 (half of the magic constant 34) is called a complementary pair. Link-lines are geometric patterns showing how these pairs are positioned within the magic square.
\end{definition}

\begin{remark}
The concept of link-lines is not new and has long been known in the study of magic squares. It is mentioned in the works of Henry Ernest Dudeney \cite{omori2018,andrews1976,benson1976} and has played an important role in understanding the structure of magic squares.
\end{remark}

In a fourth-order magic square, when dividing the numbers 1 to 16 into eight pairs, the combinations that sum to 17 are:
\begin{align*}
&(1, 16), (2, 15), (3, 14), (4, 13), \\
&(5, 12), (6, 11), (7, 10), (8, 9)
\end{align*}
The "link-line type" of a magic square is determined by how these eight pairs are spatially arranged.

\subsection{Equivalence of Link-Lines}
We consider the essential structure to be unchanged by rotations (90, 180, 270 degrees) or reflections. Therefore, link-line patterns that can be transformed into one another by these symmetric operations are treated as the same type. Specifically, link-line types are normalized by identifying orbits under the action of the dihedral group $D_4$ (order 8).

\section{Classification of 880 Magic Squares into 12 Types}

\subsection{Overview of Classification Results}
By classifying all 880 types of fourth-order magic squares by their link-line patterns, they are consolidated into \textbf{12 distinct types}. Table \ref{tab:summary} shows the frequency of each type. To visually represent the link-line patterns, we created structural diagrams where pairs summing to 17 are indicated by the same letters (A--H).

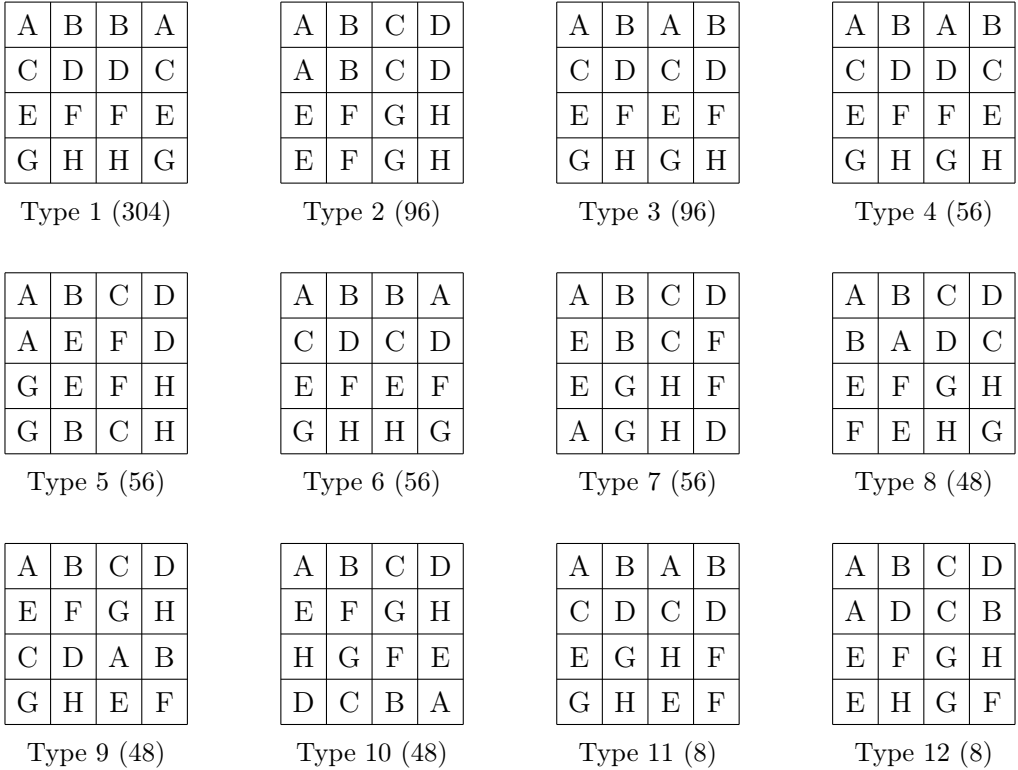
\begin{figure}[H]
\centering
% Type 1
\begin{minipage}{0.22\textwidth}
\centering
\begin{tikzpicture}[scale=0.6]
\foreach \i in {0,1,2,3,4} { \draw (0,\i) -- (4,\i); \draw (\i,0) -- (\i,4); }
\node at (0.5,3.5) {A}; \node at (1.5,3.5) {B}; \node at (2.5,3.5) {B}; \node at (3.5,3.5) {A};
\node at (0.5,2.5) {C}; \node at (1.5,2.5) {D}; \node at (2.5,2.5) {D}; \node at (3.5,2.5) {C};
\node at (0.5,1.5) {E}; \node at (1.5,1.5) {F}; \node at (2.5,1.5) {F}; \node at (3.5,1.5) {E};
\node at (0.5,0.5) {G}; \node at (1.5,0.5) {H}; \node at (2.5,0.5) {H}; \node at (3.5,0.5) {G};
\end{tikzpicture}
\\ {\small Type 1 (304)}
\end{minipage}
% Type 2
\begin{minipage}{0.22\textwidth}
\centering
\begin{tikzpicture}[scale=0.6]
\foreach \i in {0,1,2,3,4} { \draw (0,\i) -- (4,\i); \draw (\i,0) -- (\i,4); }
\node at (0.5,3.5) {A}; \node at (1.5,3.5) {B}; \node at (2.5,3.5) {C}; \node at (3.5,3.5) {D};
\node at (0.5,2.5) {A}; \node at (1.5,2.5) {B}; \node at (2.5,2.5) {C}; \node at (3.5,2.5) {D};
\node at (0.5,1.5) {E}; \node at (1.5,1.5) {F}; \node at (2.5,1.5) {G}; \node at (3.5,1.5) {H};
\node at (0.5,0.5) {E}; \node at (1.5,0.5) {F}; \node at (2.5,0.5) {G}; \node at (3.5,0.5) {H};
\end{tikzpicture}
\\ {\small Type 2 (96)}
\end{minipage}
% Type 3
\begin{minipage}{0.22\textwidth}
\centering
\begin{tikzpicture}[scale=0.6]
\foreach \i in {0,1,2,3,4} { \draw (0,\i) -- (4,\i); \draw (\i,0) -- (\i,4); }
\node at (0.5,3.5) {A}; \node at (1.5,3.5) {B}; \node at (2.5,3.5) {A}; \node at (3.5,3.5) {B};
\node at (0.5,2.5) {C}; \node at (1.5,2.5) {D}; \node at (2.5,2.5) {C}; \node at (3.5,2.5) {D};
\node at (0.5,1.5) {E}; \node at (1.5,1.5) {F}; \node at (2.5,1.5) {E}; \node at (3.5,1.5) {F};
\node at (0.5,0.5) {G}; \node at (1.5,0.5) {H}; \node at (2.5,0.5) {G}; \node at (3.5,0.5) {H};
\end{tikzpicture}
\\ {\small Type 3 (96)}
\end{minipage}
% Type 4
\begin{minipage}{0.22\textwidth}
\centering
\begin{tikzpicture}[scale=0.6]
\foreach \i in {0,1,2,3,4} { \draw (0,\i) -- (4,\i); \draw (\i,0) -- (\i,4); }
\node at (0.5,3.5) {A}; \node at (1.5,3.5) {B}; \node at (2.5,3.5) {A}; \node at (3.5,3.5) {B};
\node at (0.5,2.5) {C}; \node at (1.5,2.5) {D}; \node at (2.5,2.5) {D}; \node at (3.5,2.5) {C};
\node at (0.5,1.5) {E}; \node at (1.5,1.5) {F}; \node at (2.5,1.5) {F}; \node at (3.5,1.5) {E};
\node at (0.5,0.5) {G}; \node at (1.5,0.5) {H}; \node at (2.5,0.5) {G}; \node at (3.5,0.5) {H};
\end{tikzpicture}
\\ {\small Type 4 (56)}
\end{minipage}

\vspace{1.5em}

% Type 5
\begin{minipage}{0.22\textwidth}
\centering
\begin{tikzpicture}[scale=0.6]
\foreach \i in {0,1,2,3,4} { \draw (0,\i) -- (4,\i); \draw (\i,0) -- (\i,4); }
\node at (0.5,3.5) {A}; \node at (1.5,3.5) {B}; \node at (2.5,3.5) {C}; \node at (3.5,3.5) {D};
\node at (0.5,2.5) {A}; \node at (1.5,2.5) {E}; \node at (2.5,2.5) {F}; \node at (3.5,2.5) {D};
\node at (0.5,1.5) {G}; \node at (1.5,1.5) {E}; \node at (2.5,1.5) {F}; \node at (3.5,1.5) {H};
\node at (0.5,0.5) {G}; \node at (1.5,0.5) {B}; \node at (2.5,0.5) {C}; \node at (3.5,0.5) {H};
\end{tikzpicture}
\\ {\small Type 5 (56)}
\end{minipage}
% Type 6
\begin{minipage}{0.22\textwidth}
\centering
\begin{tikzpicture}[scale=0.6]
\foreach \i in {0,1,2,3,4} { \draw (0,\i) -- (4,\i); \draw (\i,0) -- (\i,4); }
\node at (0.5,3.5) {A}; \node at (1.5,3.5) {B}; \node at (2.5,3.5) {B}; \node at (3.5,3.5) {A};
\node at (0.5,2.5) {C}; \node at (1.5,2.5) {D}; \node at (2.5,2.5) {C}; \node at (3.5,2.5) {D};
\node at (0.5,1.5) {E}; \node at (1.5,1.5) {F}; \node at (2.5,1.5) {E}; \node at (3.5,1.5) {F};
\node at (0.5,0.5) {G}; \node at (1.5,0.5) {H}; \node at (2.5,0.5) {H}; \node at (3.5,0.5) {G};
\end{tikzpicture}
\\ {\small Type 6 (56)}
\end{minipage}
% Type 7
\begin{minipage}{0.22\textwidth}
\centering
\begin{tikzpicture}[scale=0.6]
\foreach \i in {0,1,2,3,4} { \draw (0,\i) -- (4,\i); \draw (\i,0) -- (\i,4); }
\node at (0.5,3.5) {A}; \node at (1.5,3.5) {B}; \node at (2.5,3.5) {C}; \node at (3.5,3.5) {D};
\node at (0.5,2.5) {E}; \node at (1.5,2.5) {B}; \node at (2.5,2.5) {C}; \node at (3.5,2.5) {F};
\node at (0.5,1.5) {E}; \node at (1.5,1.5) {G}; \node at (2.5,1.5) {H}; \node at (3.5,1.5) {F};
\node at (0.5,0.5) {A}; \node at (1.5,0.5) {G}; \node at (2.5,0.5) {H}; \node at (3.5,0.5) {D};
\end{tikzpicture}
\\ {\small Type 7 (56)}
\end{minipage}
% Type 8
\begin{minipage}{0.22\textwidth}
\centering
\begin{tikzpicture}[scale=0.6]
\foreach \i in {0,1,2,3,4} { \draw (0,\i) -- (4,\i); \draw (\i,0) -- (\i,4); }
\node at (0.5,3.5) {A}; \node at (1.5,3.5) {B}; \node at (2.5,3.5) {C}; \node at (3.5,3.5) {D};
\node at (0.5,2.5) {B}; \node at (1.5,2.5) {A}; \node at (2.5,2.5) {D}; \node at (3.5,2.5) {C};
\node at (0.5,1.5) {E}; \node at (1.5,1.5) {F}; \node at (2.5,1.5) {G}; \node at (3.5,1.5) {H};
\node at (0.5,0.5) {F}; \node at (1.5,0.5) {E}; \node at (2.5,0.5) {H}; \node at (3.5,0.5) {G};
\end{tikzpicture}
\\ {\small Type 8 (48)}
\end{minipage}

\vspace{1.5em}

% Type 9
\begin{minipage}{0.22\textwidth}
\centering
\begin{tikzpicture}[scale=0.6]
\foreach \i in {0,1,2,3,4} { \draw (0,\i) -- (4,\i); \draw (\i,0) -- (\i,4); }
\node at (0.5,3.5) {A}; \node at (1.5,3.5) {B}; \node at (2.5,3.5) {C}; \node at (3.5,3.5) {D};
\node at (0.5,2.5) {E}; \node at (1.5,2.5) {F}; \node at (2.5,2.5) {G}; \node at (3.5,2.5) {H};
\node at (0.5,1.5) {C}; \node at (1.5,1.5) {D}; \node at (2.5,1.5) {A}; \node at (3.5,1.5) {B};
\node at (0.5,0.5) {G}; \node at (1.5,0.5) {H}; \node at (2.5,0.5) {E}; \node at (3.5,0.5) {F};
\end{tikzpicture}
\\ {\small Type 9 (48)}
\end{minipage}
% Type 10
\begin{minipage}{0.22\textwidth}
\centering
\begin{tikzpicture}[scale=0.6]
\foreach \i in {0,1,2,3,4} { \draw (0,\i) -- (4,\i); \draw (\i,0) -- (\i,4); }
\node at (0.5,3.5) {A}; \node at (1.5,3.5) {B}; \node at (2.5,3.5) {C}; \node at (3.5,3.5) {D};
\node at (0.5,2.5) {E}; \node at (1.5,2.5) {F}; \node at (2.5,2.5) {G}; \node at (3.5,2.5) {H};
\node at (0.5,1.5) {H}; \node at (1.5,1.5) {G}; \node at (2.5,1.5) {F}; \node at (3.5,1.5) {E};
\node at (0.5,0.5) {D}; \node at (1.5,0.5) {C}; \node at (2.5,0.5) {B}; \node at (3.5,0.5) {A};
\end{tikzpicture}
\\ {\small Type 10 (48)}
\end{minipage}
% Type 11
\begin{minipage}{0.22\textwidth}
\centering
\begin{tikzpicture}[scale=0.6]
\foreach \i in {0,1,2,3,4} { \draw (0,\i) -- (4,\i); \draw (\i,0) -- (\i,4); }
\node at (0.5,3.5) {A}; \node at (1.5,3.5) {B}; \node at (2.5,3.5) {A}; \node at (3.5,3.5) {B};
\node at (0.5,2.5) {C}; \node at (1.5,2.5) {D}; \node at (2.5,2.5) {C}; \node at (3.5,2.5) {D};
\node at (0.5,1.5) {E}; \node at (1.5,1.5) {G}; \node at (2.5,1.5) {H}; \node at (3.5,1.5) {F};
\node at (0.5,0.5) {G}; \node at (1.5,0.5) {H}; \node at (2.5,0.5) {E}; \node at (3.5,0.5) {F};
\end{tikzpicture}
\\ {\small Type 11 (8)}
\end{minipage}
% Type 12
\begin{minipage}{0.22\textwidth}
\centering
\begin{tikzpicture}[scale=0.6]
\foreach \i in {0,1,2,3,4} { \draw (0,\i) -- (4,\i); \draw (\i,0) -- (\i,4); }
\node at (0.5,3.5) {A}; \node at (1.5,3.5) {B}; \node at (2.5,3.5) {C}; \node at (3.5,3.5) {D};
\node at (0.5,2.5) {A}; \node at (1.5,2.5) {D}; \node at (2.5,2.5) {C}; \node at (3.5,2.5) {B};
\node at (0.5,1.5) {E}; \node at (1.5,1.5) {F}; \node at (2.5,1.5) {G}; \node at (3.5,1.5) {H};
\node at (0.5,0.5) {E}; \node at (1.5,0.5) {H}; \node at (2.5,0.5) {G}; \node at (3.5,0.5) {F};
\end{tikzpicture}
\\ {\small Type 12 (8)}
\end{minipage}

\caption{Geometric patterns of 12 Connection Line Types}
\label{fig:types}
\end{figure}

\begin{table}[H]
\centering
\caption{Count of Magic Squares by Link-Line Type}
\label{tab:summary}
\begin{tabular}{cc}
\toprule
Type & Count \\
\midrule
1 & 304 \\
2, 3 & 96 \\
4--7 & 56 \\
8--10 & 48 \\
11, 12 & 8 \\
\bottomrule
\end{tabular}
\end{table}

A complete list of indices for fourth-order magic squares belonging to each link-line type is provided in Appendix A.

\section{Connecting Link-Line Classification to APD Analysis}
The classification by link-line types introduced in this section captures the geometric structure of fourth-order magic squares discretely and provides the foundation for the Alternating Power Difference (APD) classification treated in subsequent sections. In fact, while APD is defined as an algebraic invariant of the matrix representation of a magic square, the distribution and equivalence of its values are subject to significant structural constraints for each link-line type.

In this sense, link-line classification is a natural preprocessing step preceding APD classification, providing a framework to clarify the correspondence between geometric symmetry and algebraic properties.

\section{What is the Alternating Power Difference (APD)?}

\subsection{A New Metric for Measuring Matrix Symmetry}
To understand the internal structure of magic squares more deeply, this study introduces an algebraic invariant called the \textbf{Alternating Power Difference (APD)}. The APD is a new concept that quantifies high-order symmetries within a matrix through its interaction with the symmetric group. In this study, an algebraic invariant refers to a quantity that remains invariant under simultaneous row and column permutations $A \mapsto P_\tau A P_\tau^{-1}$ ($\tau \in S_n$).

\subsection{Definition of APD}
Let $S_n$ be the symmetric group on the set $\{1, \ldots, n\}$, and let $\operatorname{sgn}: S_n \to \{\pm 1\}$ be the sign homomorphism \cite{rotman1995,cameron1999}.

\begin{definition}[Function derived from a matrix]
For an $n \times n$ matrix $A$ and a permutation $\sigma \in S_n$, define the function $f_A: S_n \to \mathbb{Z}$ as:
\[
f_A(\sigma) := \operatorname{tr}(AP_\sigma) = \sum_{i=1}^n A_{i,\sigma(i)}
\]
where $P_\sigma$ is the permutation matrix corresponding to $\sigma$.
\end{definition}

This quantity is the sum of the elements of $A$ along the permutation $\sigma$, based on the same principle as those appearing in the definition of the permanent or determinant \cite{minc1978,brualdi2006}.

\begin{definition}[Alternating Power Difference]
\label{def:apd_magic}
For an integer-valued function $f: S_n \to \mathbb{Z}$ and $m \geq 1$, the Alternating Power Difference is defined as:
\[
\operatorname{APD}_m(f) := \sum_{\sigma \in S_n} \operatorname{sgn}(\sigma) f(\sigma)^m
\]
For a matrix $A$, we denote $\operatorname{APD}_m(A) := \operatorname{APD}_m(f_A)$.
\end{definition}

\begin{remark}
$\operatorname{APD}_m(A) = 0$ is equivalent to saying that the sums of $f_A(\sigma)^m$ over even and odd permutations are perfectly equal:
\[
\sum_{\sigma \in A_n} f_A(\sigma)^m = \sum_{\sigma \notin A_n} f_A(\sigma)^m
\]
where $A_n = \{\sigma \in S_n \mid \operatorname{sgn}(\sigma) = +1\}$ is the alternating group.
\end{remark}

\subsection{First Appearance Degree and Value}

\begin{definition}[First Appearance Degree]
For a function $f: S_n \to \mathbb{Z}$, the smallest $m \geq 1$ such that $\operatorname{APD}_m(f) \neq 0$ is called the \textbf{First Appearance Degree} of $f$, denoted by $m_1(f)$. If $\operatorname{APD}_m(f) = 0$ for all $m \geq 1$, we define $m_1(f) = \infty$.
\end{definition}

\begin{definition}[First Appearance Value]
When the first appearance degree $m_1(f)$ is finite, $\operatorname{APD}_{m_1}(f)$ is called the \textbf{First Appearance Value}.
\end{definition}

The first appearance degree $m_1$ represents the "minimum degree at which the multisets generated by even and odd permutations become distinguishable" and is an important invariant characterizing the algebraic complexity of the matrix. The multisets referred to here are $\{ f(\sigma) \mid \sigma \in A_n \}$ and $\{ f(\sigma) \mid \sigma \notin A_n \}$.

\section{Classification of Fourth-Order Magic Squares by APD}

\subsection{Classification Methodology}
We performed APD analysis on all 880 fourth-order magic squares following these steps:
\begin{enumerate}
\item Treat each magic square as a $4 \times 4$ matrix $A$.
\item For all 24 permutations $\sigma$ of the symmetric group $S_4$, calculate $f_A(\sigma) = \sum_{i=1}^4 A_{i,\sigma(i)}$.
\item Calculate $\operatorname{APD}_m(A)$ for $m = 1, 2, 3, \ldots$ and record the smallest $m$ where it is non-zero as the first appearance degree $m_1$.
\item Group magic squares with identical absolute first appearance values $|\operatorname{APD}_{m_1}(A)|$.
\item Sort the groups by $m_1$, then by ascending absolute value, and index them as G1, G2, etc.
\end{enumerate}
Since sign differences correspond to specific symmetric transformations of the matrix (such as transposition or reflection), we focus on the essential algebraic structure by using absolute values in this classification.

\subsection{Classification Results}
As a result of APD classification, all 880 types of fourth-order magic squares are divided into 51 groups based on their first appearance degree and value. Details for each group are provided in Appendix B.

\subsection{Distribution and Characteristics of First Appearance Degree}
Based on the first appearance degree $m_1$, the 880 magic squares fall into three classes:
\begin{itemize}
\item \textbf{Class $m_1 = 3$}: 240 squares
\item \textbf{Class $m_1 = 4$}: 624 squares
\item \textbf{Class $m_1 = \infty$ (G51)}: 16 squares
\end{itemize}

\subsubsection{Meaning of $m_1 = 3$ and $m_1 = 4$}
For magic squares with $m_1 = 4$, the following holds:
\[
\operatorname{APD}_1(A) = \operatorname{APD}_2(A) = \operatorname{APD}_3(A) = 0
\]
This means the contributions of even and odd permutations are perfectly balanced up to the third degree, indicating an inherent higher-order algebraic equilibrium. In contrast, for squares with $m_1 = 3$, the APD first becomes non-zero at the third degree, showing asymmetry at a lower level.

\subsubsection{G51: Complete Alternating Magic Squares}
The 16 magic squares in group G51 satisfy $\operatorname{APD}_m(A) = 0$ for all $m \geq 1$. That is, the multisets of $f_A(\sigma)$ for even and odd permutations are identical, and $m_1 = \infty$. In this study, we name these \textbf{Complete Alternating Magic Squares}. They represent magic squares in a state of purely algebraic perfect balance, undetectable by geometric classifications such as link-lines.

\subsection{Discrete Structure of APD Values}
Observation of Table \ref{tab:magic_squares} shows that APD values are not continuously distributed but exhibit a clear discrete structure, concentrating on specific integer values. Especially in the $m_1 = 4$ groups, the APD values are distributed across a finite set of discrete points, reflecting the strong constraints on component arrangement under the action of the permutation group $S_4$.

\subsection{Integer Structure in the $m_1 = 3$ Groups}
Investigating the absolute APD values for magic squares in groups with $m_1 = 3$ (G1--G11), we confirmed that all are divisible by a common factor:
\[
768 = 2^8 \times 3
\]
The values normalized by this factor are shown in Table \ref{tab:m1_3_ratios}.

\begin{table}[H]
\centering
\caption{Normalization of APD Values in $m_1 = 3$ Groups}
\label{tab:m1_3_ratios}
\begin{tabular}{cccc}
\toprule
Group & $|\operatorname{APD}|$ & $|\operatorname{APD}|/768$ & Count \\
\midrule
G1  & 1,536  & 2  & 40 \\
G2  & 2,304  & 3  & 16 \\
G3  & 3,072  & 4  & 24 \\
G4  & 3,840  & 5  & 8 \\
G5  & 4,608  & 6  & 48 \\
G6  & 5,376  & 7  & 8 \\
G7  & 6,144  & 8  & 16 \\
G8  & 6,912  & 9  & 8 \\
G9  & 7,680  & 10 & 16 \\
G10 & 9,216  & 12 & 24 \\
G11 & 12,288 & 16 & 32 \\
\bottomrule
\end{tabular}
\end{table}

The set of values obtained after normalization is $\{2,3,4,5,6,7,8,9,10,12,16\}$, with 11, 13, 14, and 15 noticeably absent. This demonstrates that $m_1 = 3$ magic squares possess an extremely restricted integer structure from the perspective of APD.

\subsection{Supplementary Observation}
This set of integers coincides with the set of possible orders for torsion subgroups of elliptic curves over the field of rational numbers, as known from Mazur's theorem \cite{mazur1977,silverman2009}. The corresponding sequence is registered in the OEIS as A059765 \cite{OEISA059765}. We record this as a numerical observation; whether this suggests a common principle between the APD structure of magic squares and algebraic constraints in other fields remains a subject for future research.

\section{Correspondence Between APD and Link-Line Classifications}
This section analyzes the relationship between the 12 geometric link-line types (L1--L12) and the 51 algebraic APD groups (G1--G51). Detailed tables are in Appendix C; here we focus on observed \textbf{structural features}.

\subsection{Coarse Separation by First Appearance Degree}
The most striking fact is that link-line types are largely separated by the first appearance degree $m_1$:
\begin{itemize}
\item Groups with $m_1 = 3$ belong primarily to L4--L7.
\item Groups with $m_1 = 4$ are distributed across L1--L3 and L8--L10.
\item The group with $m_1 = \infty$ (Complete Alternating Magic Squares) belongs only to L1.
\end{itemize}
This separation is almost without exception, showing that $m_1$ is strongly correlated with geometric classification.

\subsection{Appearance of Perfectly Uniform Distributions}
In all APD groups belonging to $m_1 = 3$, the distribution across link-line types L4, L5, L6, and L7 is always perfectly uniform (each type appears with the same frequency within each group). Similarly, in some $m_1 = 4$ groups, a perfectly uniform distribution is observed across L8, L9, and L10. This suggests that sets of magic squares with specific APD values possess high symmetry even within the link-line classification.

\subsection{Exceptional Types: L11 and L12}
Link-line types L11 and L12 are the least frequent. Interestingly, they appear only in $m_1 = 3$ groups and never in groups where $m_1 \geq 4$. This limited appearance suggests that geometrically special structures do not necessarily entail increased algebraic complexity, and that APD classification serves to refine link-line classification.

\subsection{Summary}
The analysis reveals that:
\begin{itemize}
\item $m_1$ is an effective metric strongly linked to geometric classification.
\item APD classification visualizes symmetries and balanced structures hidden within link-line classification.
\item Geometric frequency and algebraic equilibrium are not independent but partially correlated.
\end{itemize}
These correspondences show that geometric perspectives and algebraic invariants play complementary roles in classifying magic squares.

\section{Complete Alternating Magic Squares}

\subsection{Definition and Characteristics}
The 16 magic squares of group G51 possess the extremely special property that $\operatorname{APD}_m = 0$ for all degrees $m \geq 1$. In this section, we define these as \textbf{Complete Alternating Magic Squares} and analyze their algebraic structure.

\begin{figure}[H]
\centering
\begin{tabular}{cccc}
% Row 1
\begin{minipage}{0.22\textwidth}
\centering
\textbf{No.2}\\[0.5em]
\begin{tabular}{|c|c|c|c|} \hline 1 & 2 & 15 & 16 \\ \hline 13 & 14 & 3 & 4 \\ \hline 12 & 7 & 10 & 5 \\ \hline 8 & 11 & 6 & 9 \\ \hline \end{tabular}
\end{minipage} &
\begin{minipage}{0.22\textwidth}
\centering
\textbf{No.87}\\[0.5em]
\begin{tabular}{|c|c|c|c|} \hline 1 & 7 & 14 & 12 \\ \hline 9 & 15 & 4 & 6 \\ \hline 8 & 2 & 13 & 11 \\ \hline 16 & 10 & 3 & 5 \\ \hline \end{tabular}
\end{minipage} &
\begin{minipage}{0.22\textwidth}
\centering
\textbf{No.100}\\[0.5em]
\begin{tabular}{|c|c|c|c|} \hline 1 & 8 & 9 & 16 \\ \hline 14 & 13 & 4 & 3 \\ \hline 7 & 2 & 15 & 10 \\ \hline 12 & 11 & 6 & 5 \\ \hline \end{tabular}
\end{minipage} &
\begin{minipage}{0.22\textwidth}
\centering
\textbf{No.179}\\[0.5em]
\begin{tabular}{|c|c|c|c|} \hline 1 & 12 & 13 & 8 \\ \hline 15 & 10 & 3 & 6 \\ \hline 2 & 7 & 14 & 11 \\ \hline 16 & 5 & 4 & 9 \\ \hline \end{tabular}
\end{minipage} \\[4em]

% Row 2
\begin{minipage}{0.22\textwidth}
\centering
\textbf{No.211}\\[0.5em]
\begin{tabular}{|c|c|c|c|} \hline 2 & 1 & 16 & 15 \\ \hline 14 & 13 & 4 & 3 \\ \hline 11 & 8 & 9 & 6 \\ \hline 7 & 12 & 5 & 10 \\ \hline \end{tabular}
\end{minipage} &
\begin{minipage}{0.22\textwidth}
\centering
\textbf{No.285}\\[0.5em]
\begin{tabular}{|c|c|c|c|} \hline 2 & 7 & 10 & 15 \\ \hline 11 & 12 & 5 & 6 \\ \hline 8 & 1 & 16 & 9 \\ \hline 13 & 14 & 3 & 4 \\ \hline \end{tabular}
\end{minipage} &
\begin{minipage}{0.22\textwidth}
\centering
\textbf{No.318}\\[0.5em]
\begin{tabular}{|c|c|c|c|} \hline 2 & 8 & 11 & 13 \\ \hline 10 & 16 & 5 & 3 \\ \hline 7 & 1 & 12 & 14 \\ \hline 15 & 9 & 6 & 4 \\ \hline \end{tabular}
\end{minipage} &
\begin{minipage}{0.22\textwidth}
\centering
\textbf{No.376}\\[0.5em]
\begin{tabular}{|c|c|c|c|} \hline 2 & 11 & 14 & 7 \\ \hline 16 & 9 & 4 & 5 \\ \hline 1 & 8 & 13 & 12 \\ \hline 15 & 6 & 3 & 10 \\ \hline \end{tabular}
\end{minipage} \\[4em]

% Row 3
\begin{minipage}{0.22\textwidth}
\centering
\textbf{No.448}\\[0.5em]
\begin{tabular}{|c|c|c|c|} \hline 3 & 4 & 13 & 14 \\ \hline 15 & 16 & 1 & 2 \\ \hline 6 & 9 & 8 & 11 \\ \hline 10 & 5 & 12 & 7 \\ \hline \end{tabular}
\end{minipage} &
\begin{minipage}{0.22\textwidth}
\centering
\textbf{No.482}\\[0.5em]
\begin{tabular}{|c|c|c|c|} \hline 3 & 6 & 15 & 10 \\ \hline 13 & 8 & 1 & 12 \\ \hline 4 & 9 & 16 & 5 \\ \hline 14 & 11 & 2 & 7 \\ \hline \end{tabular}
\end{minipage} &
\begin{minipage}{0.22\textwidth}
\centering
\textbf{No.617}\\[0.5em]
\begin{tabular}{|c|c|c|c|} \hline 4 & 3 & 14 & 13 \\ \hline 15 & 10 & 7 & 2 \\ \hline 6 & 5 & 12 & 11 \\ \hline 9 & 16 & 1 & 8 \\ \hline \end{tabular}
\end{minipage} &
\begin{minipage}{0.22\textwidth}
\centering
\textbf{No.619}\\[0.5em]
\begin{tabular}{|c|c|c|c|} \hline 4 & 3 & 14 & 13 \\ \hline 16 & 15 & 2 & 1 \\ \hline 5 & 10 & 7 & 12 \\ \hline 9 & 6 & 11 & 8 \\ \hline \end{tabular}
\end{minipage} \\[4em]

% Row 4
\begin{minipage}{0.22\textwidth}
\centering
\textbf{No.645}\\[0.5em]
\begin{tabular}{|c|c|c|c|} \hline 4 & 5 & 16 & 9 \\ \hline 14 & 7 & 2 & 11 \\ \hline 3 & 10 & 15 & 6 \\ \hline 13 & 12 & 1 & 8 \\ \hline \end{tabular}
\end{minipage} &
\begin{minipage}{0.22\textwidth}
\centering
\textbf{No.665}\\[0.5em]
\begin{tabular}{|c|c|c|c|} \hline 4 & 6 & 15 & 9 \\ \hline 14 & 12 & 7 & 1 \\ \hline 3 & 5 & 10 & 16 \\ \hline 13 & 11 & 2 & 8 \\ \hline \end{tabular}
\end{minipage} &
\begin{minipage}{0.22\textwidth}
\centering
\textbf{No.777}\\[0.5em]
\begin{tabular}{|c|c|c|c|} \hline 5 & 3 & 16 & 10 \\ \hline 11 & 13 & 8 & 2 \\ \hline 6 & 4 & 9 & 15 \\ \hline 12 & 14 & 1 & 7 \\ \hline \end{tabular}
\end{minipage} &
\begin{minipage}{0.22\textwidth}
\centering
\textbf{No.793}\\[0.5em]
\begin{tabular}{|c|c|c|c|} \hline 5 & 6 & 11 & 12 \\ \hline 16 & 9 & 8 & 1 \\ \hline 3 & 4 & 13 & 14 \\ \hline 10 & 15 & 2 & 7 \\ \hline \end{tabular}
\end{minipage}
\end{tabular}
\vspace{1.5em}
\caption{Complete Alternating Magic Squares 1--16 (G51 Group)}
\label{fig:perfect_alternating}
\end{figure}

\begin{definition}[Complete Alternating Magic Square]
A fourth-order magic square $A$ is a \textbf{Complete Alternating Magic Square} if, for all positive integers $m \geq 1$:
\[
\operatorname{APD}_m(A) = \sum_{\sigma \in S_4} \operatorname{sgn}(\sigma) f_A(\sigma)^m = 0
\]
\end{definition}

\begin{observation}
For a Complete Alternating Magic Square $A$, the condition that $\operatorname{APD}_m(A) = 0$ for all $m$ is equivalent to: the multiset generated by even permutations $\{f_A(\sigma) \mid \sigma \in A_4\}$ and the multiset generated by odd permutations $\{f_A(\sigma) \mid \sigma \notin A_4\}$ are identical.
\end{observation}

This observation shows that in a Complete Alternating Magic Square, the 12 even and 12 odd permutations possess a \textbf{perfectly symmetric structure} regarding the values of $f_A$.

\section{Analysis of Specific Examples of Complete Alternating Magic Squares}
This section analyzes the distribution of $f_A(\sigma)$ values for magic square No. 2 as a representative example.

\subsection{Representative Case: Complete Alternating Magic Square No. 2}

\begin{table}[H]
\centering
\caption{Complete Alternating Magic Square No. 2}
\begin{tabular}{|c|c|c|c|}
\hline 1 & 2 & 15 & 16 \\ \hline 13 & 14 & 3 & 4 \\ \hline 12 & 7 & 10 & 5 \\ \hline 8 & 11 & 6 & 9 \\ \hline
\end{tabular}
\end{table}

This square satisfies $\operatorname{APD}_m(A)=0$ for all $m \ge 1$.

\subsection{Distribution of $f_A(\sigma)$ Values}
Calculating $f_A(\sigma)$ for all 24 permutations of $S_4$, we find 9 distinct values with the following frequencies:

\begin{table}[H]
\centering
\caption{Values and Frequencies of $f_A(\sigma)$ (No. 2)}
\label{tab:f_distribution_no2}
\begin{tabular}{ccc}
\toprule
$f$ value & Even Permutations & Odd Permutations \\
\midrule
18 & 1 & 1 \\
20 & 1 & 1 \\
24 & 1 & 1 \\
26 & 2 & 2 \\
34 & 2 & 2 \\
42 & 2 & 2 \\
44 & 1 & 1 \\
48 & 1 & 1 \\
50 & 1 & 1 \\
\midrule
Total & 12 & 12 \\
\bottomrule
\end{tabular}
\end{table}

The frequencies for even and odd permutations match perfectly for every value.

\subsection{Common Distribution Structures}
Numerical verification shows that the 16 Complete Alternating Magic Squares fall into \textbf{two types} regarding the distribution of $f_A(\sigma)$:
\begin{itemize}
\item \textbf{Type I}: No. 2, 179, 211, 376, 448, 482, 619, 645
\item \textbf{Type II}: No. 87, 100, 285, 318, 617, 665, 777, 793
\end{itemize}
While the sets of $f$ values differ between types, the distribution patterns (multiplicities) are identical.

\section{Classification via Symmetries in Matrix Products}
Beyond APD, these squares exhibit symmetry in matrix products. We investigated three types: $A^2$, $A \cdot \mathrm{rot}_{180}(A)$, and $A \cdot A^T$.

\subsection{Rotationally Symmetric Type: Equality of Self-Product and Rotation-Product}
Exactly \textbf{8 squares} satisfy:
\begin{equation}
A^2 = A \cdot \mathrm{rot}_{180}(A)
\end{equation}
This indicates that the point symmetry of the component arrangement extends to the matrix's multiplicative structure (Table \ref{tab:rotation_symmetric}).

\begin{table}[H]
\centering
\caption{Squares where Self-Product equals Rotation-Product}
\label{tab:rotation_symmetric}
\begin{tabular}{ccc}
\toprule
Index & $\mathrm{tr}(A^2)$ & $\mathrm{tr}(A A^T)$ \\
\midrule
2   & 1236 & 1496 \\
100 & 1356 & 1496 \\
211 & 1140 & 1496 \\
285 & 1356 & 1496 \\
448 & 1236 & 1496 \\
617 & 972  & 1496 \\
619 & 1140 & 1496 \\
793 & 972  & 1496 \\
\bottomrule
\end{tabular}
\end{table}

\subsection{Dual Type: Pair Structure Sharing Transpose Products}

The remaining 8 squares form four dual pairs $(A_1, A_2)$ that satisfy the following identity:
\begin{equation}
A_1 A_1^T = A_2 A_2^T
\end{equation}
The specific pairs are listed in Table \ref{tab:transpose_pairs}. While the transpose product $AA^T$ is always a symmetric matrix \cite{horn2012}, it is extremely rare for distinct magic squares to yield identical transpose products. This suggests a deep algebraic correspondence between the squares in each pair.

\begin{table}[H]
\centering
\caption{Dual pairs sharing identical transpose products}
\label{tab:transpose_pairs}
\begin{tabular}{ccc}
\toprule
Pair & Index & Index \\
\midrule
1 & 87  & 777 \\
2 & 179 & 645 \\
3 & 318 & 665 \\
4 & 376 & 482 \\
\bottomrule
\end{tabular}
\end{table}

\subsection{Universal Invariants: The Trace}

Notably, for \textbf{all 16 Complete Alternating Magic Squares}, the trace of the transpose product remains invariant, taking a constant value:
\begin{equation}
\mathrm{tr}(AA^T) = 1496
\end{equation}
In contrast, the trace of the matrix square, $\mathrm{tr}(A^2)$, takes one of four distinct values, with each appearing exactly four times, as shown in Table \ref{tab:trace_distribution}.

\begin{table}[H]
\centering
\caption{Distribution of the trace of the matrix square $\mathrm{tr}(A^2)$}
\label{tab:trace_distribution}
\begin{tabular}{ccc}
\toprule
$\mathrm{tr}(A^2)$ Value & Frequency & Magic Square IDs \\
\midrule
972  & 4 & 617, 665, 777, 793 \\
1140 & 4 & 211, 376, 619, 645 \\
1236 & 4 & 2, 179, 448, 482 \\
1356 & 4 & 87, 100, 285, 318 \\
\bottomrule
\end{tabular}
\end{table}

This perfectly uniform distribution suggests the existence of structural invariance that cannot be explained by mere coincidence.

\subsection{Double Bifurcation Structure}

Based on the results above, the Complete Alternating Magic Squares are partitioned into two distinct categories from the perspective of matrix products:
\begin{itemize}
    \item \textbf{Rotational Symmetry Type (8 squares)}: Satisfying $A^2 = A \cdot \mathrm{rot}_{180}(A)$.
    \item \textbf{Dual Type (8 squares)}: Forming dual pairs that share identical transpose products.
\end{itemize}
It is important to emphasize that this classification is independent of the bifurcation based on $f_A(\sigma)$ discussed in the previous section, differing in both definition and origin. In other words, the Complete Alternating Magic Squares possess a \textbf{multi-layered structure characterized by both combinatorial and matrix-algebraic symmetries}.

\section{Algebraic Descent in Squared Matrices}

Complete Alternating Magic Squares possess an extremely strong algebraic equilibrium, satisfying $\operatorname{APD}_m(A)=0$ for all degrees $m \ge 1$. In this section, we systematically investigate the APD behavior of the squared matrix
\[
A^2 = A \cdot A
\]
for magic squares $A$ possessing this property.

The focus of this section is to elucidate how "infinite-order perfect alternation" transforms through the matrix power operation—specifically, to reveal the reality of the \textbf{algebraic descent from perfect equilibrium}.

\subsection{Investigation Setup}

The subjects of this study are the 16 Complete Alternating Magic Squares $A$ extracted from the 880 magic squares of order 4. For each $A$, we calculated the APD of $A^2$ for each degree and recorded the first occurrence degree
\[
m_1(A^2)
\]
at which the APD becomes non-zero, along with its corresponding value.

\subsection{Uniform First Occurrence Degree}

As a result of numerical computation, it was confirmed that the following holds for all 16 Complete Alternating Magic Squares:

\begin{observation}[Algebraic Descent in Squaring]
For a Complete Alternating Magic Square $A$, its squared matrix $A^2$ possesses the following properties:
\begin{enumerate}
    \item $m_1(A^2)=4$
    \item $\operatorname{APD}_1(A^2) = \operatorname{APD}_2(A^2) = \operatorname{APD}_3(A^2)=0$
    \item $\operatorname{APD}_4(A^2)\neq 0$
\end{enumerate}
\end{observation}

In other words, the Complete Alternating Magic Squares, which originally had $m_1(A)=\infty$, uniformly undergo a descent in their first occurrence degree through the squaring operation:
\[
\infty \longrightarrow 4
\]
This behavior is common across all 16 squares, with no exceptions.

\subsection{Distribution Structure of APD$_4$ Values}

Next, examining the values of $\operatorname{APD}_4(A^2)$, we find that the 16 magic squares are classified into \textbf{four distinct values} (Table \ref{tab:apd4_distribution}).

\begin{table}[H]
\centering
\caption{Distribution of $\operatorname{APD}_4$ values in squared matrices}
\label{tab:apd4_distribution}
\begin{tabular}{ccc}
\toprule
$\operatorname{APD}_4(A^2)$ value & Frequency & Corresponding Magic Squares \\
\midrule
$+230,400,000$ & 4 & No. 87, 100, 285, 318 \\
$-16,515,072$  & 4 & No. 211, 376, 619, 645 \\
$-294,912,000$ & 4 & No. 2, 179, 448, 482 \\
$-429,023,232$ & 4 & No. 617, 665, 777, 793 \\
\bottomrule
\end{tabular}
\end{table}

Three points are particularly noteworthy here:
\begin{enumerate}
    \item \textbf{Perfectly Uniform Distribution}: The four types of values appear exactly four times each.
    \item \textbf{Discreteness of Values}: The $\text{APD}_4$ values are not continuous but are strictly limited to these four values.
    \item \textbf{Sign Structure}: There is only one type of positive value, while the remaining three are negative.
\end{enumerate}

\subsection{Correspondence with Matrix Product Symmetry}

Comparing the distribution of $\operatorname{APD}_4$ values with the \textbf{bifurcation based on matrix products} (Rotational Symmetry Type and Dual Type) revealed in the previous section, the following facts are confirmed:
\begin{itemize}
    \item In both groups, the four types of $\operatorname{APD}_4$ values appear \textbf{exactly twice each}.
    \item Therefore, the distribution of $\operatorname{APD}_4$ values is consistent with the bifurcation by matrix products.
\end{itemize}
This implies that classifications derived from \textbf{different principles} than the bifurcation based on $f_A(\sigma)$ ultimately converge toward the same numerical equilibrium structure.

\subsection{Discussion}

From the results above, we can conclude the following:
\begin{enumerate}
    \item Complete Alternating Magic Squares always descend to $m_1=4$ through the squaring operation.
    \item This descent is not a variation between individual magic squares but a phenomenon common to the entire class.
    \item The $\text{APD}_4$ values possess a strong equilibrium structure consisting of 4 types with 4 instances each.
    \item The bifurcation by matrix product symmetry and the classification by APD reflect the same structural constraints from different perspectives.
\end{enumerate}

In this sense, the Complete Alternating Magic Squares occupy the highest level of a clear hierarchical structure with "perfect equilibrium" at its apex. Their squared matrices are understood as objects that have algebraically descended while partially retaining that structure.

\section{Future Research Topics}

This study has revealed that by using the algebraic invariant known as the Alternating Power Difference (APD), a singular class called \emph{Complete Alternating Magic Squares} exists within the set of all order-4 magic squares—a class not captured by conventional classifications. Simultaneously, this discovery presents the following essential unsolved problems:

\begin{enumerate}
    \item \textbf{Existence in Higher Orders}: Do Complete Alternating Magic Squares exist for $n \ge 5$? In particular, determining whether order 4 is a unique case is of great importance.
    \item \textbf{Establishment of Construction Theory}: Is there a method to systematically construct the 16 Complete Alternating Magic Squares of order 4 from symmetry or algebraic conditions?
    \item \textbf{General Behavior under Power Operations}: For a Complete Alternating Magic Square $A$, what is the general behavior of the first occurrence degree $m_1(A^k)$ for $A^k$?
    \item \textbf{Structural Meaning of APD Values}: Can the discrete numerical values and factor structures appearing in APD values be explained group-theoretically or combinatorially?
    \item \textbf{Extension to Other Fields}: Can the concept of APD be extended to Latin squares or other combinatorial arrangements?
\end{enumerate}

These are all fundamental tasks for establishing APD not merely as a computational quantity but as a structural theory.

\section{Conclusion}

In this study, we introduced a new algebraic invariant called the Alternating Power Difference (APD) for all 880 types of order-4 magic squares and performed a systematic classification based on the first occurrence degree $m_1$.

As a result, the following points were clarified:
\begin{enumerate}
    \item Classification by APD provides a new perspective for capturing the algebraic structure of magic squares, independent of conventional geometric or visual classifications.
    \item There exist 16 \emph{Complete Alternating Magic Squares} with $m_1 = \infty$, forming a special class with extremal parity symmetry.
    \item Complete Alternating Magic Squares exhibit a strong equilibrium structure even in matrix products and powers, allowing them to be understood as algebraic objects that transcend mere numerical arrangements.
\end{enumerate}

In conclusion, APD introduces a new framework connecting matrix theory, group theory, and combinatorics to the study of magic squares. The discovery of Complete Alternating Magic Squares is the first concrete example demonstrating the theoretical potential of this framework.

We hope that this study will deepen the understanding of symmetry in magic squares and serve as a starting point for the study of more general combinatorial structures.

\section*{Data Availability}

The data for the 880 types of $4 \times 4$ magic squares used in this study is provided as a CSV file named \texttt{magic\_squares\_4x4.csv}. 
In this file, each magic square is assigned a serial number (ID) from 1 to 880, and the 16 integer components of each square are recorded in row-major order. 
These IDs correspond one-to-one with the classification results based on APD (Alternating Power Difference), group numbers, and the identification of perfectly alternating magic squares discussed in this paper.

Furthermore, all primary numerical calculations and verifications used in this study were executed using Python programs. 
To ensure transparency and reproducibility, the CSV file and all source code used for the analysis are publicly available in the following GitHub repository:
\begin{quote}
    \url{https://github.com/soaisu-ken/magic-square-apd-classification/releases/tag/v1.0.0}
\end{quote}

Calculations based on the APD definition formula, the determination of the first occurrence degree $m_1$, the classification of all 880 squares, and the verification of the distribution results claimed in this paper are fully reproducible using the provided code, enabling independent verification by third parties.

\appendix
\clearpage
\section*{Appendix A:List of Magic Square IDs by Connection Line Type}

The following is an exhaustive list of magic square IDs belonging to each connection line type.

\subsubsection*{Type 1 (304 squares)}
\begin{small}
1, 2, 4, 5, 6, 7, 12, 13, 14, 15, 16, 17, 18, 19, 36, 37, 38, 39, 44, 45, 46, 47, 48, 49, 50, 51, 52, 53, 54, 55, 70, 71, 72, 73, 74, 75, 76, 77, 78, 79, 80, 81, 86, 87, 91, 100, 110, 111, 119, 121, 123, 125, 127, 130, 131, 132, 133, 134, 135, 136, 137, 138, 139, 143, 144, 150, 151, 152, 153, 154, 155, 156, 157, 158, 166, 167, 168, 179, 182, 184, 186, 187, 188, 189, 190, 194, 198, 199, 210, 211, 218, 219, 220, 221, 222, 223, 224, 225, 226, 227, 252, 253, 254, 255, 256, 257, 258, 259, 261, 262, 284, 285, 286, 287, 288, 295, 296, 298, 307, 309, 310, 311, 312, 313, 314, 315, 318, 320, 336, 337, 338, 339, 341, 342, 352, 359, 363, 366, 367, 370, 376, 378, 379, 380, 381, 382, 383, 391, 394, 397, 398, 404, 427, 428, 429, 430, 432, 433, 435, 436, 447, 448, 452, 454, 456, 457, 458, 459, 474, 475, 477, 480, 482, 488, 490, 498, 499, 502, 516, 520, 521, 533, 534, 540, 543, 545, 547, 556, 557, 561, 573, 574, 593, 594, 595, 596, 597, 599, 600, 602, 603, 604, 609, 610, 611, 612, 614, 615, 617, 619, 624, 625, 626, 627, 629, 630, 631, 633, 634, 636, 638, 639, 641, 645, 651, 653, 654, 656, 665, 667, 670, 671, 672, 673, 680, 683, 694, 696, 697, 699, 706, 707, 708, 716, 718, 719, 725, 739, 740, 742, 743, 745, 747, 750, 751, 752, 755, 756, 760, 761, 763, 764, 769, 771, 772, 774, 776, 777, 780, 791, 792, 793, 801, 802, 804, 806, 807, 813, 815, 816, 819, 821, 822, 824, 825, 826, 830, 832, 833, 836, 837, 840, 841, 842, 846, 848, 849, 851, 856, 857, 859, 861, 862, 864
\end{small}

\subsubsection*{Type 2 (96 squares)}
\begin{small}
24, 25, 30, 31, 60, 61, 64, 65, 84, 85, 92, 93, 140, 141, 146, 147, 159, 160, 162, 163, 191, 192, 195, 196, 216, 217, 236, 237, 250, 251, 273, 274, 321, 322, 325, 326, 333, 334, 350, 351, 384, 385, 386, 387, 399, 400, 401, 402, 419, 420, 443, 444, 451, 466, 467, 504, 507, 508, 515, 527, 528, 548, 549, 550, 566, 569, 570, 572, 581, 582, 589, 590, 649, 663, 664, 669, 681, 682, 705, 712, 713, 715, 722, 723, 730, 731, 734, 735, 766, 778, 809, 810, 817, 843, 854, 855
\end{small}

\subsubsection*{Type 3 (96 squares)}
\begin{small}
32, 33, 34, 35, 66, 67, 68, 69, 95, 96, 97, 98, 101, 103, 106, 108, 114, 115, 148, 149, 164, 165, 169, 173, 228, 229, 230, 231, 264, 265, 266, 267, 280, 282, 291, 293, 301, 302, 328, 329, 330, 331, 344, 345, 356, 362, 388, 389, 423, 424, 425, 426, 460, 461, 462, 463, 468, 471, 479, 481, 511, 512, 513, 514, 523, 524, 553, 554, 575, 576, 577, 578, 620, 622, 640, 642, 657, 658, 659, 660, 674, 675, 676, 677, 709, 710, 720, 721, 784, 787, 797, 798, 827, 838, 852, 853
\end{small}

\subsubsection*{Type 4 (56 squares)}
\begin{small}
8, 9, 10, 11, 20, 26, 40, 41, 42, 43, 58, 94, 99, 128, 129, 170, 212, 232, 242, 243, 244, 245, 248, 275, 276, 277, 278, 327, 332, 353, 354, 357, 409, 410, 411, 412, 417, 437, 440, 495, 496, 497, 541, 542, 579, 585, 586, 605, 606, 607, 608, 686, 687, 688, 726, 727
\end{small}

\subsubsection*{Type 5 (56 squares)}
\begin{small}
23, 29, 59, 105, 118, 142, 145, 161, 172, 180, 193, 197, 200, 205, 207, 208, 215, 235, 249, 272, 283, 335, 348, 349, 358, 373, 390, 403, 405, 406, 407, 408, 418, 441, 442, 509, 526, 551, 555, 559, 563, 567, 568, 571, 580, 587, 588, 644, 701, 711, 728, 729, 732, 733, 736, 737
\end{small}

\subsubsection*{Type 6 (56 squares)}
\begin{small}
238, 239, 240, 241, 260, 263, 340, 343, 413, 414, 415, 416, 431, 434, 491, 492, 493, 494, 500, 501, 519, 522, 544, 546, 598, 601, 616, 618, 652, 655, 684, 685, 753, 754, 759, 762, 770, 773, 794, 796, 800, 811, 812, 814, 820, 823, 847, 863, 865, 866, 871, 872, 874, 876, 877, 880
\end{small}

\subsubsection*{Type 7 (56 squares)}
\begin{small}
268, 271, 300, 303, 346, 347, 369, 371, 438, 439, 453, 455, 470, 472, 486, 510, 517, 518, 525, 529, 531, 538, 552, 564, 692, 693, 700, 703, 714, 717, 738, 749, 757, 758, 765, 767, 775, 781, 782, 783, 786, 795, 799, 805, 829, 831, 845, 858, 867, 868, 869, 870, 873, 875, 878, 879
\end{small}

\subsubsection*{Type 8 (48 squares)}
\begin{small}
21, 22, 27, 28, 56, 57, 62, 63, 82, 83, 89, 90, 213, 214, 233, 234, 246, 247, 269, 270, 316, 317, 323, 324, 421, 422, 445, 446, 450, 464, 465, 503, 505, 506, 583, 584, 591, 592, 648, 661, 662, 668, 678, 679, 768, 779, 818, 844
\end{small}

\subsubsection*{Type 9 (48 squares)}
\begin{small}
102, 104, 107, 109, 116, 117, 171, 174, 177, 178, 201, 204, 279, 281, 292, 294, 304, 305, 355, 365, 372, 375, 393, 396, 469, 473, 483, 485, 530, 532, 536, 537, 560, 565, 621, 623, 646, 647, 690, 691, 702, 704, 744, 748, 785, 788, 828, 839
\end{small}

\subsubsection*{Type 10 (48 squares)}
\begin{small}
112, 113, 120, 122, 124, 126, 175, 176, 183, 185, 203, 206, 289, 290, 297, 299, 306, 308, 360, 361, 368, 377, 392, 395, 476, 478, 487, 489, 535, 539, 558, 562, 628, 632, 635, 637, 695, 698, 741, 746, 789, 790, 803, 808, 834, 835, 850, 860
\end{small}

\subsubsection*{Type 11 (8 squares)}
\begin{small}
3, 88, 209, 319, 449, 613, 650, 666
\end{small}

\subsubsection*{Type 12 (8 squares)}
\begin{small}
181, 202, 364, 374, 484, 643, 689, 724
\end{small}

\section*{Appendix B:Classification of Order-4 Magic Squares into 51 Groups by APD}

Table \ref{tab:magic_squares} shows the results of the classification into 51 groups using APD. For each group, the first occurrence degree $m_1$, the absolute value of the first occurrence $|\operatorname{APD}|$, and the number of magic squares included are listed.

\begin{longtable}{lccc}
\caption{Classification of 880 types of $4 \times 4$ magic squares into 51 groups} \label{tab:magic_squares} \\
\toprule
Group & $m_1$ & $|\operatorname{APD}|$ & Count \\
\midrule
\endfirsthead

\multicolumn{4}{c}%
{{\bfseries \tablename\ \thetable{} -- Continued}} \\
\toprule
Group & $m_1$ & $|\operatorname{APD}|$ & Count \\
\midrule
\endhead

\midrule
\multicolumn{4}{r}{{Continued on next page}} \\
\endfoot

\bottomrule
\endlastfoot

G1 & 3 & 1,536 & 40 \\
G2 & 3 & 2,304 & 16 \\
G3 & 3 & 3,072 & 24 \\
G4 & 3 & 3,840 & 8 \\
G5 & 3 & 4,608 & 48 \\
G6 & 3 & 5,376 & 8 \\
G7 & 3 & 6,144 & 16 \\
G8 & 3 & 6,912 & 8 \\
G9 & 3 & 7,680 & 16 \\
G10 & 3 & 9,216 & 24 \\
G11 & 3 & 12,288 & 32 \\
G12 & 4 & 2,112 & 4 \\
G13 & 4 & 20,928 & 4 \\
G14 & 4 & 29,568 & 4 \\
G15 & 4 & 38,400 & 16 \\
G16 & 4 & 43,008 & 24 \\
G17 & 4 & 47,616 & 24 \\
G18 & 4 & 49,152 & 24 \\
G19 & 4 & 50,112 & 4 \\
G20 & 4 & 55,296 & 48 \\
G21 & 4 & 57,408 & 4 \\
G22 & 4 & 58,368 & 4 \\
G23 & 4 & 59,520 & 16 \\
G24 & 4 & 64,512 & 4 \\
G25 & 4 & 65,088 & 4 \\
G26 & 4 & 79,872 & 48 \\
G27 & 4 & 91,392 & 16 \\
G28 & 4 & 97,920 & 4 \\
G29 & 4 & 102,912 & 24 \\
G30 & 4 & 103,488 & 4 \\
G31 & 4 & 104,448 & 16 \\
G32 & 4 & 141,312 & 4 \\
G33 & 4 & 147,840 & 4 \\
G34 & 4 & 156,672 & 4 \\
G35 & 4 & 168,960 & 24 \\
G36 & 4 & 183,744 & 4 \\
G37 & 4 & 237,120 & 4 \\
G38 & 4 & 251,328 & 4 \\
G39 & 4 & 267,840 & 4 \\
G40 & 4 & 288,768 & 24 \\
G41 & 4 & 290,304 & 48 \\
G42 & 4 & 296,448 & 24 \\
G43 & 4 & 314,496 & 4 \\
G44 & 4 & 344,064 & 24 \\
G45 & 4 & 345,600 & 48 \\
G46 & 4 & 370,176 & 48 \\
G47 & 4 & 388,608 & 24 \\
G48 & 4 & 430,080 & 24 \\
G49 & 4 & 448,512 & 4 \\
G50 & 4 & 457,728 & 4 \\
G51 & $\infty$ & 0 & 16 \\
\midrule
\multicolumn{3}{l}{\textbf{Total}} & \textbf{880} \\
\end{longtable}

\section*{Appendix C:Correspondence between APD Classification and Connection Line Classification}

The 880 types of order-4 magic squares were classified from two different perspectives: classification into 12 types based on geometric connection lines (L1--L12) and classification into 51 groups based on algebraic APD (G1--G51). 
We investigated which connection line type (L1--L12) the magic squares belonging to each APD group (G1--G51) fall into. By clarifying this correspondence, it is possible to understand the interaction between geometric structure and algebraic invariants.

\subsection{Overview of the Correspondence}

Table \ref{tab:apd_connection_distribution} shows the distribution of connection line types within each APD group. From this table, the following remarkable structures can be observed.

\begin{longtable}{ccp{8cm}}
\caption{Correspondence between APD groups and connection line types} \label{tab:apd_connection_distribution} \\
\toprule
APD Group & Total Count & Connection Line Type Distribution \\
\midrule
\endfirsthead

\multicolumn{3}{c}%
{{\bfseries \tablename\ \thetable{} -- Continued}} \\
\toprule
APD Group & Total Count & Connection Line Type Distribution \\
\midrule
\endhead

\midrule
\multicolumn{3}{r}{{Continued on next page}} \\
\endfoot

\bottomrule
\endlastfoot

\multicolumn{3}{l}{\textbf{Groups with First Occurrence Degree $m_1 = 3$}} \\
\midrule
G1  & 40 & L4:10, L5:10, L6:10, L7:10 \\
G2  & 16 & L4:4, L5:4, L6:4, L7:4 \\
G3  & 24 & L4:6, L5:6, L6:6, L7:6 \\
G4  & 8  & L4:2, L5:2, L6:2, L7:2 \\
G5  & 48 & L4:10, L5:10, L6:10, L7:10, \textcolor{red}{L11:4, L12:4} \\
G6  & 8  & L4:2, L5:2, L6:2, L7:2 \\
G7  & 16 & L4:4, L5:4, L6:4, L7:4 \\
G8  & 8  & L4:2, L5:2, L6:2, L7:2 \\
G9  & 16 & L4:4, L5:4, L6:4, L7:4 \\
G10 & 24 & L4:4, L5:4, L6:4, L7:4, \textcolor{red}{L11:4, L12:4} \\
G11 & 32 & L4:8, L5:8, L6:8, L7:8 \\
\midrule
\multicolumn{3}{l}{\textbf{Groups with First Occurrence Degree $m_1 = 4$}} \\
\midrule
G12 & 4  & L1:4 \\
G13 & 4  & L1:4 \\
G14 & 4  & L1:4 \\
G15 & 16 & L1:16 \\
G16 & 24 & L1:8, L2:8, L3:8 \\
G17 & 24 & L1:8, L2:8, L3:8 \\
G18 & 24 & L1:8, L2:8, L3:8 \\
G19 & 4  & L1:4 \\
G20 & 48 & \textcolor{blue}{L8:16, L9:16, L10:16} \\
G21 & 4  & L1:4 \\
G22 & 4  & L1:4 \\
G23 & 16 & L1:16 \\
G24 & 4  & L1:4 \\
G25 & 4  & L1:4 \\
G26 & 48 & L1:16, L2:16, L3:16 \\
G27 & 16 & L1:16 \\
G28 & 4  & L1:4 \\
G29 & 24 & L1:8, L2:8, L3:8 \\
G30 & 4  & L1:4 \\
G31 & 16 & L1:16 \\
G32 & 4  & L1:4 \\
G33 & 4  & L1:4 \\
G34 & 4  & L1:4 \\
G35 & 24 & L1:24 \\
G36 & 4  & L1:4 \\
G37 & 4  & L1:4 \\
G38 & 4  & L1:4 \\
G39 & 4  & L1:4 \\
G40 & 24 & L1:8, L2:8, L3:8 \\
G41 & 48 & \textcolor{blue}{L8:16, L9:16, L10:16} \\
G42 & 24 & L1:8, L2:8, L3:8 \\
G43 & 4  & L1:4 \\
G44 & 24 & L1:8, L2:8, L3:8 \\
G45 & 48 & \textcolor{blue}{L8:16, L9:16, L10:16} \\
G46 & 48 & L1:16, L2:16, L3:16 \\
G47 & 24 & L1:8, L2:8, L3:8 \\
G48 & 24 & L1:24 \\
G49 & 4  & L1:4 \\
G50 & 4  & L1:4 \\
\midrule
\multicolumn{3}{l}{\textbf{Groups with First Occurrence Degree $m_1 = \infty$}} \\
\midrule
G51 & 16 & \textcolor{red}{L1:16} \\
\end{longtable}

\end{document}